# STEPUP PROCEDURES CONTROLLING GENERALIZED FWER AND GENERALIZED FDR[1]

### By Sanat K. Sarkar

#### *Temple University*

In many applications of multiple hypothesis testing where more than one false rejection can be tolerated, procedures controlling error rates measuring at least $k$ false rejections, instead of at least one, for some fixed $k \geq 1$ can potentially increase the ability of a procedure to detect false null hypotheses. The $k$-FWER, a generalized version of the usual familywise error rate (FWER), is such an error rate that has recently been introduced in the literature and procedures controlling it have been proposed. A further generalization of a result on the $k$-FWER is provided in this article. In addition, an alternative and less conservative notion of error rate, the $k$-FDR, is introduced in the same spirit as the $k$-FWER by generalizing the usual false discovery rate (FDR). A $k$-FWER procedure is constructed given any set of increasing constants by utilizing the $k$th order joint null distributions of the $p$-values without assuming any specific form of dependence among all the $p$-values. Procedures controlling the $k$-FDR are also developed by using the $k$th order joint null distributions of the $p$-values, first assuming that the sets of null and nonnull $p$-values are mutually independent or they are jointly positively dependent in the sense of being multivariate totally positive of order two (MTP$_2$) and then discarding that assumption about the overall dependence among the $p$-values.

**1. Introduction.** Having realized that the traditional idea of controlling the familywise error rate (FWER), which is the probability of rejecting at least one true null hypothesis, is too stringent to use when a large number of hypotheses are simultaneously tested, researchers have focused in the last decade on defining alternative less stringent error rates and developing methods that control them. The false discovery rate (FDR), which is the

Received January 2006; revised January 2007.

[1]Supported by NSF Grants DMS-03-06366 and DMS-06-03868.

expected proportion of falsely rejected null hypotheses and which was introduced by Benjamini and Hochberg [1], is the first of these that has received considerable attention [2, 3, 4, 5, 6, 12, 15, 16, 17, 21, 22]. Recently, the ideas of controlling the probabilities of falsely rejecting at least $k$ null hypotheses, which is the $k$-FWER, and the false discovery proportion (FDP) exceeding a certain threshold $\gamma \in [0, 1)$ have been introduced as alternatives to the FWER and methods controlling these new error rates have been suggested [10, 11, 13, 18, 23].

Sarkar [18] developed single-step and stepwise $k$-FWER procedures utilizing the $k$th order joint null distributions of the test statistics. He generalized the Bonferroni single-step procedure and obtained its Holm [8] type improvement, thereby providing a generalized $k$-FWER stepdown procedure. These are different from the corresponding procedures of Lehmann and Romano [11]. He then generalized the procedure of Hochberg [7] by using a stepup procedure with the same critical values as those of his generalized Holm procedure. He proved that this generalized version of Hochberg's procedure controls the $k$-FWER by making use of a generalized Simes' inequality he obtained for statistics that are positively dependent in the sense of being multivariate totally positive of order two (MTP$_2$), a condition due to Karlin and Rinott [9] and often shared by test statistics in multiple testing. These alternative procedures are often more powerful than those in [11] based on marginal distributions, especially when the statistics are close to being independent.

In this article we continue our research in the line of [18] and develop newer procedures using $k$th order joint null distributions of the test statistics. First, we generalize the method of Romano and Shaikh [13] and construct a stepup $k$-FWER procedure given any set of increasing constants. Second, which is more interesting, we introduce a less conservative notion of error rate than the $k$-FWER, which is the $k$-FDR (to be defined later in this section), by generalizing the usual FDR in the same spirit as the $k$-FWER, and provide newer stepup procedures that control it. The $k$-FWER procedure is constructed without assuming any specific overall dependence among the test statistics. The $k$-FDR procedures are derived under two scenarios–first, under the assumption that either the test statistics corresponding to the true null hypotheses are independent of those corresponding to the false null hypotheses or they are jointly MTP$_2$, and then without that assumption.

Let us denote by $V$ and $R$ the total number of false rejections and the total number of rejections, respectively, of null hypotheses. Then

$$(1.1) \qquad\qquad k\text{-FWER} = Pr\{V \geq k\}.$$

We will generalize it further in this article in terms of the following measure:

$$(1.2) \quad k\text{-FDR} = E(k\text{-FDP}), \quad \text{where } k\text{-FDP} = \begin{cases} \dfrac{V}{R}, & \text{if } V \geq k, \\ 0, & \text{otherwise.} \end{cases}$$

The concept of $k$-FDR has not been considered before, as far as we know, even though consideration of it, instead of the $k$-FWER, appears to be a more natural extension of the idea of using the FDR as a less restrictive error rate than the FWER. Clearly, it reduces to the usual FDR when $k = 1$, and, more importantly, as $k$-FDR $\leq$ $k$-FWER, controlling it would be a less conservative approach than controlling the $k$-FWER.

The construction of the $k$-FWER stepup procedure is provided in Section 2. Section 3 is devoted to the derivation of the $k$-FDR stepup procedures. The $k$-FDR procedures are generalized versions of the usual FDR procedures of Benjamini and Hochberg [1] and Benjamini and Yekutieli [2], referred to as the generalized BH and generalized BY procedures, respectively, in this paper. The generalized BH version of the $k$-FDR procedure provides uniformly better control of the $k$-FDR, as it is intended to do so, than the generalized Hochberg procedure in [18] that, being a $k$-FWER procedure, controls the $k$-FDR under the same distributional setup. As $k$-FDR $\leq$ FDR, an FDR procedure, such as the original BH procedure, also controls the $k$-FDR. However, as our simulation studies indicate, although the generalized BH $k$-FDR procedure does not seem to outperform the original BH procedure when the number of true null hypotheses is small, its performance is much better when this number is relatively large and the test statistics are not highly dependent on each other.

Sarkar [18] generalized Simes' test [20] by controlling the probability of at least $k$, instead of one, false rejections under the intersection of the null hypotheses. Interestingly, unlike the usual BH procedure, its generalized version that controls the $k$-FDR is not based on these generalized Simes' critical values. In fact, we prove that the stepup procedure with the generalized Simes' critical values does not control the $k$-FDR.

Before we proceed to develop procedures with a control of the $k$-FWER or $k$-FDR, we recall here the definitions of stepdown and stepup procedures. Consider testing $n$ null hypotheses $H_1, \ldots, H_n$ simultaneously against certain alternatives using their $p$-values $P_1, \ldots, P_n$, respectively. Let $P_{1:n} \leq \cdots \leq P_{n:n}$ denote the ordered $p$-values. Then, given some critical values $\alpha_1 \leq \cdots \leq \alpha_n$, a stepdown procedure accepts $H_i$ for all $i \geq j_{SD}$ and rejects the rest, where $j_{SD} = \min_{1 \leq i \leq n}\{i : P_{i:n} \geq \alpha_i\}$, if the minimum exists; otherwise, it rejects all the hypotheses. A stepup procedure, on the other hand, rejects $H_i$ for all $i \leq j_{SU}$ and accepts the rest, where $j_{SU} = \max_{1 \leq i \leq n}\{i : P_{i:n} \leq \alpha_i\}$, if the maximum exists; otherwise, it will accept all the hypotheses. These can be generalized by considering $\alpha_1 = \cdots = \alpha_k$, for some fixed $1 \leq k \leq n$.

**2. $k$-FWER controlling stepup procedure.** In this section we consider developing a stepup procedure with a control of the $k$-FWER at $\alpha$ starting with any increasing set of constants and using the $k$th order joint null distribution of the $p$-values. This is an attempt to generalize the idea of Romano and Shaikh [13] that uses only the marginal $p$-values. First, we recall the following inequality from [18] that holds for any set of random variables $X_1, \ldots, X_n$, not necessarily $p$-values.

LEMMA 2.1.  *Let* $\mathcal{C}_k = \{J : J \subseteq \{1, \ldots, n\}, \ |J| = k\}$ *and* $a_i = \binom{i}{k}$, $i = k, \ldots, n$. *Then, given any set of constants* $c_k \leq \cdots \leq c_n$, *and* $1 \leq k \leq n$, *we have*

$$
\begin{aligned}
Pr\left\{\bigcup_{i=k}^{n}(X_{i:n} \leq c_i)\right\} &\leq \sum_{J \in \mathcal{C}_k} Pr\left\{\max_{j \in J} X_j \leq c_k\right\} \\
&\quad + \sum_{i=k+1}^{n} a_i^{-1} \sum_{J \in \mathcal{C}_k} Pr\left\{c_{i-1} < \max_{j \in J} X_j \leq c_i\right\}.
\end{aligned}
\tag{2.1}
$$

REMARK 2.1.  The above lemma generalizes Lemma 3.1 of Lehmann and Romano [11]. When the $k$th-order joint distributions are identical with $G_k(x) = Pr\{\max_{i \in J} X_i \leq x\}$ for all $J \in \mathcal{C}_k$, then it reduces to

$$
Pr\left\{\bigcup_{i=k}^{n}(X_{i:n} \leq c_i)\right\} \leq \binom{n}{k}\left[G_k(c_k) + \sum_{i=k+1}^{n} a_i^{-1}\{G_k(c_i) - G_k(c_{i-1})\}\right].
\tag{2.2}
$$

Considering $k = 1$ and uniform $p$-values, one gets the inequality given in [11].

We are now ready to describe our method of constructing a $k$-FWER stepup procedure given any set of increasing constants $\alpha_1' \leq \cdots \leq \alpha_n'$. We will, however, assume that the $k$th-order joint null distributions of the $p$-values are identical. Let $F_k(x) = Pr\{\max_{i \in J} P_i \leq x\}$ for all $J \in \mathcal{C}_k$ and

$$
D_{k,n}' = \max_{k \leq n_0 \leq n} S_{k,n}'(n_0),
\tag{2.3}
$$

where

$$
\begin{aligned}
S_{k,n}'(n_0) &= \binom{n_0}{k}\left[F_k(\alpha_{n-n_0+k}') \right. \\
&\quad \left. + \sum_{i=k+1}^{n_0} a_i^{-1}\{F_k(\alpha_{n-n_0+i}') - F_k(\alpha_{n-n_0+i-1}')\}\right]
\end{aligned}
\tag{2.4}
$$

and the probabilities are determined under the null hypotheses.

THEOREM 2.1.  *Given any set of constants* $\alpha_k' \leq \cdots \leq \alpha_n'$ *and* $0 < \alpha < 1$, *consider the stepup procedure with the critical values* $\alpha_1 \leq \cdots \leq \alpha_n$ *satisfying* $F_k(\alpha_i) = \alpha F_k(\alpha_{i \vee k}')/D_{k,n}'$, $i = 1, \ldots, n$, *where* $i \vee k = \max(i, k)$. *This controls the* $k$-FWER *at* $\alpha$.

Proof. Assume without any loss of generality that the first $n_0$ $p$-values correspond to the true null hypotheses. Let $P_{1:n_0} \leq \cdots \leq P_{n_0:n_0}$ be the ordered versions of these $p$-values. Assume that $n_0 \geq k$; otherwise, the $k$-FWER is zero and, hence, is trivially controlled. From [13] and using inequality (2.2), we see that the $k$-FWER of the stepup procedure in the theorem satisfies

$$k\text{-FWER} \leq Pr\left\{\bigcup_{i=k}^{n_0}(P_{i:n_0} \leq \alpha_{n-n_0+i})\right\}$$

$$(2.5) \quad \leq \binom{n_0}{k}\left[F_k(\alpha_{n-n_0+k}) + \sum_{i=k+1}^{n_0} a_i^{-1}\{F_k(\alpha_{n-n_0+i}) - F_k(\alpha_{n-n_0+i-1})\}\right]$$

$$= \alpha \frac{S'_{k,n}(n_0)}{D'_{k,n}} \leq \alpha,$$

which proves the theorem. □

Remark 2.2. Sarkar [18] proposed a stepdown procedure with the critical values $\alpha_1 \leq \cdots \leq \alpha_n$ satisfying $F_k(\alpha_i) = \alpha/a_{n+k-i\vee k}$, $i = 1, \ldots, n$, and proved that it controls the $k$-FWER without assuming any overall dependence structure among the test statistics. It is alternative to, and often more powerful than, the stepdown procedure with the critical values $\alpha_i = k\alpha/(n + k - i \vee k)$, $i = 1, \ldots, n$, that Lehmann and Romano [11] proposed. Sarkar [18] also showed that a stepup procedure with the same critical values as those of his $k$-FWER stepdown procedure can also control the $k$-FWER when the test statistics are all jointly $\text{MTP}_2$. This generalizes Hochberg's procedure and its FWER control property under similar positive dependence condition; see [7, 14, 19]. When such $\text{MTP}_2$ condition does not hold for the overall joint distribution of the test statistics, Theorem 2.1 says that one can appropriately rescale these critical values before using the stepup procedure in order to control the $k$-FWER.

## 3. $k$-FDR controlling stepup procedures.

In this section we will construct stepup procedures that control the $k$-FDR using the $k$th-order joint null distribution of the $p$-values under two different scenarios, first when the sets of null and non-null $p$-values are known to be either mutually independent or are jointly dependent in the sense being $\text{MTP}_2$ and second when no such condition is known.

### 3.1. Generalized BH procedure.

Toward developing the generalized BH version of our $k$-FDR procedure, we first obtain a result providing a convenient expression for the upper bound of the $k$-FDR of a general stepup procedure. To that end, we will be using the following notation.

For a stepup procedure with the critical values $\alpha_1 \leq \cdots \leq \alpha_n$, let $R_{n-q}^{(-i_1, \ldots, -i_q)}$ and $V_{n-q}^{(-i_1, \ldots, -i_q)}$ denote respectively the number of rejections and the number of false rejections of null hypotheses when the stepup procedure based on the subset of $p$-values $\{P_1, \ldots, P_n\} \setminus \{P_{i_1}, \ldots, P_{i_q}\}$ and the critical values $\alpha_{q+1} \leq \cdots \leq \alpha_n$ is used. Let $\{H_i, i \in I_0\}$ be the set of true null hypotheses ($|I_0| = n_0$). It is assumed that $n_0 \geq k$; otherwise, the $k$-FDR $= 0$, and hence, there is nothing to prove.

Moreover, we will be using the following inequality.

LEMMA 3.1. *In testing $n$ null hypotheses, of which $n_0 \geq k$ are true, let $V_n$ and $R_n$ denote the number of false rejections and the total number of rejections, respectively. Then, we have*

$$(3.1) \qquad I(R_n = r, V_n \geq k) \leq \frac{(n - r + k)V_n}{n_0 k} I(R_n = r, V_n \geq k),$$

*for all $1 \leq k \leq r \leq n$.*

PROOF. The lemma follows by combining the following two inequalities: (i) $V_n(n_0 - V_n + k) \geq n_0 k$, which holds when $k \leq V_n \leq n_0$, and (ii) $n_0 - V_n \leq n - R_n$, which is true as $R_n - V_n \leq n - n_0$. □

Since $V = \sum_{i \in I_0} I(H_i \text{ rejected})$, we have

$k$-FDP

$$= \frac{V}{R} I(V \geq k) = \sum_{r=k}^{n} \frac{1}{r} \sum_{i \in I_0} I(P_i \leq \alpha_r,\ R = r,\ V \geq k)$$

$$= \sum_{r=k}^{n} \frac{1}{r} \sum_{i \in I_0} I(P_i \leq \alpha_r,\ R_{n-1}^{(-i)} = r - 1,\ V_{n-1}^{(-i)} \geq k - 1)$$

$$\leq \sum_{r=k}^{n} \frac{(n - r + k - 1)}{r(k-1)(n_0 - 1)}$$
$$\qquad \times \sum_{i \in I_0} V_{n-1}^{(-i)} I(P_i \leq \alpha_r,\ R_{n-1}^{(-i)} = r - 1,\ V_{n-1}^{(-i)} \geq k - 1)$$

$$= \sum_{r=k}^{n} \frac{(n - r + k - 1)}{r(k-1)(n_0 - 1)}$$
$$(3.2) \qquad \times \sum_{i \neq j \in I_0} I(\max(P_i, P_j) \leq \alpha_r,\ R_{n-2}^{(-i, -j)} = r - 2,\ V_{n-2}^{(-i, -j)} \geq k - 2)$$

$$\leq$$

$$\vdots$$

$$= \sum_{r=k}^{n} \frac{(n-r+k-1)(n-r+k-2)\cdots(n-r+1)}{r(k-1)(k-2)\cdots 1(n_0-1)\cdots(n_0-k+1)}$$

$$\times \sum_{i_1 \neq \cdots \neq i_k \in I_0} I(\max(P_{i_1},\ldots,P_{i_k}) \leq \alpha_r, R_{n-k}^{(-i_1,\ldots,-i_k)} = r-k,$$

$$V_{n-k}^{(-i_1,\ldots,-i_k)} \geq 0)$$

$$= \frac{kn_0}{a_{n_0}} \sum_{r=k}^{n} \frac{a_{n+k-r}}{r(n+k-r)} \sum_{J \in \mathcal{C}_k^0} I\left(\max_{i \in J} P_i \leq \alpha_r, \ R_{n-k}^{(-J)} = r-k\right),$$

where $\mathcal{C}_k^0 = \{J : J \subseteq I_0, |J| = k\}$ and $R_{n-k}^{(-J)}$ is the number of rejections for the stepup procedure based on $\{P_i, \ i \in J^c\}$ and the critical values $\alpha_{k+1} \leq \cdots \leq \alpha_n$.

LEMMA 3.2. *For a stepup procedure with critical values $\alpha_1 \leq \cdots \leq \alpha_n$, we have*

$$k\text{-FDR} \leq \frac{n_0 a_n}{na_{n_0}} \sum_{J \in \mathcal{C}_k^0} Pr\left\{\max_{i \in J} P_i \leq \alpha_k\right\}$$

$$+ \frac{kn_0}{a_{n_0}} \sum_{r=k+1}^{n} \sum_{J \in \mathcal{C}_k^0} E\left[Pr\{R_{n-k}^{(-J)} \geq r-k \ |P_i, i \in J\}\right.$$

$$(3.3) \qquad \times \left\{\frac{a_{n+k-r}I(\max_{i \in J} P_i \leq \alpha_r)}{r(n+k-r)}\right.$$

$$\left.\left. - \frac{a_{n+k-r+1}I(\max_{i \in J} P_i \leq \alpha_{r-1})}{(r-1)(n+k-r+1)}\right\}\right].$$

PROOF. The inequality (3.2) yields

$$k\text{-FDP}$$

$$\leq \frac{kn_0}{a_{n_0}} \sum_{r=k}^{n} \frac{a_{n+k-r}}{r(n+k-r)} \sum_{J \in \mathcal{C}_k^0} I\left(\max_{i \in J} P_i \leq \alpha_r, \ R_{n-k}^{(-J)} \geq r-k\right)$$

$$- \frac{kn_0}{a_{n_0}} \sum_{r=k}^{n-1} \frac{a_{n+k-r}}{r(n+k-r)} \sum_{J \in \mathcal{C}_k^0} I\left(\max_{i \in J} P_i \leq \alpha_r, \ R_{n-k}^{(-J)} \geq r-k+1\right)$$

$$(3.4) \qquad = \frac{n_0 a_n}{na_{n_0}} \sum_{J \in \mathcal{C}_k^0} I\left(\max_{i \in J} P_i \leq \alpha_k\right)$$

$$+ \frac{kn_0}{a_{n_0}} \sum_{r=k+1}^{n} \sum_{J \in \mathcal{C}_k^0} \Bigg[ I(R_{n-k}^{(-J)} \geq r - k)$$

$$\times \Bigg\{ \frac{a_{n+k-r} I(\max_{i \in J} P_i \leq \alpha_r)}{r(n+k-r)}$$

$$- \frac{a_{n+k-r+1} I(\max_{i \in J} P_i \leq \alpha_{r-1})}{(r-1)(n+k-r+1)} \Bigg\} \Bigg],$$

from which we get the lemma by taking expectations on both sides. □

We now construct a stepup procedure that controls the $k$-FDR assuming that the $p$-values have identical $k$th-order joint null distributions and that the sets of $p$-values corresponding to the true and false null hypotheses are either independent or all the $p$-values are jointly MTP$_2$. The MTP$_2$ property is a positive dependence property that holds for many multivariate distributions arising in multiple testing, for example, multivariate normal with a common nonnegative correlation, certain mixtures of independent random variables, and so on; see, for example, [14, 15, 16, 19]. For the definition of MTP$_2$ and some of the related results to be used to prove the next result, one can see [9, 18].

THEOREM 3.1. *Consider the stepup procedure with the critical values* $\alpha_1 \leq \cdots \leq \alpha_n$ *satisfying*

$$(3.5) \qquad F_k(\alpha_i) = \begin{cases} \dfrac{\alpha}{a_n}, & \text{for } i = 1, \ldots, k, \\ \dfrac{i(n+k-i)\alpha}{kna_{n+k-i}}, & \text{for } i = k, \ldots, n. \end{cases}$$

*The $k$-FDR of this procedure is controlled at $\alpha$ if either the null $p$-values are independent of the nonnull $p$-values or the $p$-values are jointly MTP$_2$.*

PROOF. Assuming first that the null and nonnull $p$-values are mutually independent, we have from Lemma 3.2 that

$$k\text{-FDR} \leq \frac{n_0 a_n}{n} F_k(\alpha_k)$$

$$+ \frac{kn_0}{a_{n_0}} \sum_{r=k+1}^{n} \sum_{J \in \mathcal{C}_k^0} Pr\{R_{n-k}^{(-J)} \geq r - k\}$$

$$(3.6) \qquad \times \Bigg\{ \frac{a_{n+k-r} F_k(\alpha_r)}{r(n+k-r)} - \frac{a_{n+k-r+1} F_k(\alpha_{r-1})}{(r-1)(n+k-r+1)} \Bigg\},$$

which is less than or equal to $n_0 \alpha/n$, and hence controlled at $\alpha$, if $F_k(\alpha_r)$ is chosen as in (3.5) for $r = 1, \ldots, n$.

When the null and nonnull $p$-values are not mutually independent, but are jointly $\text{MTP}_2$, we will prove the theorem as follows.

Let $g(p_i, i \in J)$ be the density of $P_i, i \in J$, for any fixed $J \in \mathcal{C}_k^0$. Then, each expectation in (3.3) can be expressed as

$$(3.7)\quad
\begin{aligned}
E\bigg[ &Pr\{R_{n-k}^{(-J)} \geq r - k | P_i, i \in J\} \bigg\{ \frac{a_{n+k-r} I(\max_{i \in J} P_i \leq \alpha_r)}{r(n+k-r)} \\
&\qquad\qquad - \frac{a_{n+k-r+1} I(\max_{i \in J} P_i \leq \alpha_{r-1})}{(r-1)(n+k-r+1)} \bigg\} \bigg] \\
&= E\bigg\{ \phi_{r,J}(P_i, i \in J)\psi_{r,J}(P_i, i \in J) I\Big( \max_{i \in J} P_i \leq \alpha_r \Big) \bigg\} \\
&= F_k(\alpha_r) E^*\{\phi_{r,J}(P_i, i \in J)\psi_{r,J}(P_i, i \in J)\},
\end{aligned}$$

where the last expectation is taken with respect to the density

$$(3.8)\qquad \frac{g(p_i, i \in J) I(\max_{i \in J} p_i \leq \alpha_r)}{F_k(\alpha_r)},$$

and

$$\phi_{r,J}(p_i, i \in J) = Pr\{R_{n-k}^{(-J)} \geq r - k | P_i = p_i, i \in J\},$$

$$(3.9)\quad \psi_{r,J}(p_i, i \in J) = \frac{a_{n+k-r}}{r(n+k-r)} - \frac{a_{n+k-r+1} I(\max_{i \in J} p_i \leq \alpha_{r-1})}{(r-1)(n+k-r+1)}.$$

Since

$$(3.10)\qquad \{R_{n-k}^{(-J)} \geq r - k\} = \bigcup_{i=r-k}^{n-k} \{p_{i:J^c} \leq \alpha_{k+i}\},$$

where $p_{1:J^c} \leq \cdots \leq p_{n-k:J^c}$ are the ordered values of the set $\{p_i, i \in J^c\}$ that are all increasing in each $p_i, i \in J^c$, the indicator function $I(R_{n-k}^{(-J)} \geq r - k)$ is a decreasing (coordinatewise) function of the $p_i$'s. Also, the function $\psi_{r,J}$ is an increasing function of the $p_i$'s.

The marginal density $g(p_i, i \in J)$ of $\{P_i, i \in J\}$ is $\text{MTP}_2$. As $I(\max_{i \in J} p_i \leq \alpha_r)$ is also $\text{MTP}_2$, the density (3.8), being the product of two $\text{MTP}_2$ functions, is $\text{MTP}_2$. We will now invoke the following property of $\text{MTP}_2$ random variables. Random variables that are $\text{MTP}_2$ are positively associated; that is, any pair of functions of these random variables, both increasing or decreasing, are positively correlated. As $I(R_{n-k}^{(-J)} \geq r - k)$ is a decreasing function, the conditional probability of this function given $P_i = p_i, i \in J$, which is $\phi_{r,J}$, is a decreasing function of each $p_i, i \in J$, because of the $\text{MTP}_2$ property of the $p$-values. Therefore, we have

$$(3.11)\quad
\begin{aligned}
&E^*\{\phi_{r,J}(P_i, i \in J)\psi_{r,J}(P_i, i \in J)\} \\
&\qquad \leq E^*\{\phi_{r,J}(P_i, i \in J)\} E^*\{\psi_{r,J}(P_i, i \in J)\},
\end{aligned}$$

where

$$(3.12) \qquad \begin{aligned} &E^*\{\psi_{r,J}(P_i, i \in J)\} \\ &\quad = \frac{a_{n+k-r}}{r(n+k-r)} - \frac{a_{n+k-r+1}}{(r-1)(n+k-r+1)} \frac{F_k(\alpha_{r-1})}{F_k(\alpha_r)}, \end{aligned}$$

yielding the following inequality, again from Lemma 3.2:

$$(3.13) \qquad \begin{aligned} k\text{-FDR} &\leq \frac{n_0 a_n}{n} F_k(\alpha_k) \\ &\quad + \frac{k n_0}{a_{n_0}} \sum_{r=k+1}^{n} \sum_{J \in \mathcal{C}_k^0} E^*\{\phi_{r,J}(P_i, i \in J)\} \\ &\qquad \times \left\{ \frac{a_{n+k-r} F_k(\alpha_r)}{r(n+k-r)} - \frac{a_{n+k-r+1} F_k(\alpha_{r-1})}{(r-1)(n+k-r+1)} \right\}, \end{aligned}$$

which is less than or equal to $n_0\alpha/n$ for $F_k(\alpha_r)$ satisfying (3.5). This proves the theorem.  □

When $k = 1$, Theorem 3.1 reduces to the result establishing the usual FDR controlling property of the BH procedure; see [2, 15]. Of course, the result for the BH procedure is slightly stronger in that its FDR is less than or equal to $n_0\alpha/n$ under positive regression dependence, a slightly weaker condition than the MTP$_2$ condition, and is exactly $n_0\alpha/n$ under the independence case.

More explicit expressions of the critical values in Theorem 3.1 are

$$(3.14) \quad F_k(\alpha_i) = \begin{cases} \dfrac{k(k-1)\cdots 1\alpha}{n(n-1)\cdots(n-k+1)}, \\ \qquad \text{for } i = 1, \ldots, k, \\ \dfrac{i(k-1)(k-2)\cdots 1\alpha}{n(n-i+k-1)(n-i+k-2)\cdots(n-i+1)}, \\ \qquad \text{for } i = k, \ldots, n. \end{cases}$$

With $k = 2$,

$$(3.15) \qquad F_2(\alpha_i) = \begin{cases} \dfrac{2\alpha}{n(n-1)}, & \text{for } i = 1, 2, \\ \dfrac{i\alpha}{n(n-i+1)}, & \text{for } i = 2, \ldots, n. \end{cases}$$

In our procedure, the critical values corresponding to the smallest $k-1$ $p$-values could be chosen arbitrarily. In other words, we can always reject the null hypotheses corresponding to the smallest $k-1$ $p$-values and still control the $k$-FDR. Nevertheless, we have considered rejecting these null hypotheses based on certain critical values. By choosing these critical values all equal to $\alpha_k$, we not only preserve the monotonicity of the critical values, but also have the least conservative choice among such critical values.

It is important to see that our procedure provides uniformly better control of the $k$-FDR, as it is intended to do so, than the generalized Hochberg procedure in [18] (see Remark 2.2) that, being a $k$-FWER procedure, also controls the $k$-FDR. This is because

$$(3.16) \qquad \frac{i(n+k-i)\alpha}{nka_{n+k-i}} \geq \frac{\alpha}{a_{n+k-i}},$$

for all $i = k, \ldots, n$. The extent of this improvement can be seen in Figures 1 and 2 that are based on a numerical study to be discussed later in this section.

As mentioned in the Introduction, the original BH FDR procedure also controls the $k$-FDR. To see how it compares with our procedure, let us consider independent $p$-values and $k = 2$. The $i$th critical value of our procedure, $[i\alpha/n(n - i + 1)]^{1/2}$, exceeds the corresponding critical value $i\alpha/n$ of the BH procedure as long as $n/i(n - i + 1) \geq \alpha$, which holds for each $i$ if $4n/(n + 1)^2 \geq \alpha$. For instance, when $\alpha = 0.05$ and $n$ does not exceed 80, or when $\alpha = 0.01$ and $n$ does not exceed 400, our procedure has larger critical values, and hence, provides uniformly better control of the 2-FDR, than the BH procedure. Thus, in this case, our procedure often performs better than the BH procedure. Figures 1 and 2 provide some insight into this comparison in other cases.

Sarkar [18] has generalized Simes' test for testing the overall hypothesis $H_0 : \bigcap_{i=1}^{n} H_i$. He proposed rejecting $H_0$ if $P_{i:n} \leq \alpha_{i \vee k}$ for at least one $i = 1, \ldots, n$, for a fixed $1 \leq k \leq n$, where

$$(3.17) \quad F_k(\alpha_i) = \frac{a_i}{a_n}\alpha = \frac{i(i-1)\cdots(i-k+1)}{n(n-1)\cdots(n-k+1)}\alpha, \qquad i = k, \ldots, n,$$

and proved that it controls the probability of at least $k$ false rejections under the intersection null hypothesis at $\alpha$ exactly under independence and conservatively if the $p$-values are MTP$_2$. Interestingly, unlike what is known when $k = 1$, the generalized Simes' critical values do not always control the $k$-FDR. We provide a proof of this in the following.

Let $n_1 = n - n_0$ be the number of false null hypotheses. Consider a situation where these $n_1$ hypotheses are false to the extent that the corresponding $p$-values are all $< \alpha_1$. In other words, consider a procedure that with probability one rejects $n_1$ null hypotheses before proceeding as a stepup procedure based on the $n_0$ $p$-values that are all known to correspond to the null hypotheses and using the critical values $\alpha_{n_1+1}, \ldots, \alpha_n$. Let $R_0$ be the number of true null hypotheses that are rejected. Then the $k$-FDR of this procedure

is

$$k\text{-FDR} = E\left\{\frac{R_0}{n_1 + R_0}I(R_0 \geq k)\right\} \geq \frac{k}{n_1 + k}P\{R_0 \geq k\},$$

(3.18)

$$= \frac{k}{n_1 + k}Pr\left\{\bigcup_{i=k}^{n_0}(P_{i:I_0} \leq \alpha_{n_1+i})\right\} \geq \frac{k}{n_1 + k}Pr\{P_{k:I_0} \leq \alpha_{n_1+k}\}.$$

Let us now consider $k = 2$ and assume that the null $p$-values are i.i.d. $U(0,1)$. Then, the right-hand side of (3.18) simplifies to

$$(3.19) \qquad \frac{2}{n_1 + 2}[1 - (1 - \alpha_{n_1+2})^{n_0} - n_0\alpha_{n_1+2}(1 - \alpha_{n_1+2})^{n_0-1}],$$

which may exceed $\alpha$ for the generalized Simes critical values in some instances. For example, with $\alpha = 0.05$, $n_0 = 100$ and $n_1 = 1$, we have

$$\alpha_{n_1+2} = \left[\frac{(n_1 + 2)(n_1 + 1)\alpha}{n(n-1)}\right]^{1/2} = 0.00545,$$

from which we see that (3.19) is equal to 0.0692.

3.2. *Generalized BY procedure.* We now derive a $k$-FDR procedure only under the assumption of identical $k$th-order joint null distributions of the $p$-values, which generalizes the BY procedure in [2]. We proceed as in the construction of the generalized BH procedure, but modify those results to facilitate the development of a $k$-FDR procedure under the present scenario. First, we have the following upper bound for the $k$-FDP that follows from (3.2):

$$(3.20) \qquad k\text{-FDP} \leq \frac{kn_0a_n}{na_{n_0}}\sum_{r=k}^{n}\frac{1}{r}\sum_{J \in \mathcal{C}_k^0}I\left(\max_{i \in J}P_i \leq \alpha_r, R_{n-k}^{(-J)} = r - k\right),$$

which yields the following modification of Lemma 3.2.

LEMMA 3.3. *For a stepup procedure with critical values $\alpha_1 \leq \cdots \leq \alpha_n$, we have*

$$k\text{-FDR} \leq \frac{n_0a_n}{na_{n_0}}\sum_{J \in \mathcal{C}_k^0}Pr\left\{\max_{i \in J}P_i \leq \alpha_k\right\}$$

$$+ \frac{kn_0a_n}{na_{n_0}}\sum_{r=k+1}^{n}\sum_{J \in \mathcal{C}_k^0}E\left[Pr\{R_{n-k}^{(-J)} \geq r - k|P_i, i \in J\}\right.$$

(3.21)

$$\times\left\{\frac{I(\max_{i \in J}P_i \leq \alpha_r)}{r}\right.$$

$$\left.\left.- \frac{I(\max_{i \in J}P_i \leq \alpha_{r-1})}{r-1}\right\}\right].$$

Assuming now that the $k$-th order joint null distributions are identical, we have the following theorem.

THEOREM 3.2. *The stepup procedure with the critical values $\alpha_1 \leq \cdots \leq \alpha_n$ satisfying*

$$(3.22) \qquad F_k(\alpha_i) = \frac{(i \vee k)\alpha}{k\binom{n}{k}\sum_{r=k}^{n} 1/r}, \qquad i = 1, \ldots, n,$$

*controls the $k$-FDR at $\alpha$.*

PROOF. Lemma 3.3 yields the inequality

$$(3.23) \quad k\text{-FDR} \leq \frac{n_0 a_n}{n} F_k(\alpha_k) + \frac{k n_0 a_n}{n} \sum_{r=k+1}^{n} \frac{1}{r} [F_k(\alpha_r) - F_k(\alpha_{r-1})],$$

which is equal to $n_0 \alpha / n$ when the critical values satisfying (3.22) are used. $\square$

When $k = 1$, Theorem 3.2 reduces to the known result establishing the FDR control of the BY procedure [2]. When $k = 2$, the critical values of the generalized BY procedure are given by

$$(3.24) \qquad F_2(\alpha_i) = \frac{(i \vee 2)\alpha}{n(n-1)\sum_{r=2}^{n} 1/r}, \qquad i = 1, \ldots, n.$$

3.3. *A numerical study.* We conducted a numerical study to investigate the extent of improvement offered by the generalized BH version of our $k$-FDR procedure in controlling the $k$-FDR over the generalized Hochberg and the original BH procedures when $k = 2$.

We generated $n = 100$ dependent random variables $X_i \sim N(\mu_i, 1)$, $i = 1, \ldots, 100$, with the same variance 1 and a common correlation $\rho$, which are known to be MTP$_2$, and performed 100 hypothesis tests of $\mu = 0$ against $\mu = 2$, using each of these three procedures with $\alpha = 0.05$. The value of 2-FDP was then calculated for each procedure by setting $n_0$ of the $\mu_i$'s to zero and the remaining $n_1$ of the $\mu_i$'s to the value 2. The 2-FDR then was estimated by averaging the 2-FDP values over 5000 iterations. We did these calculations in two different scenarios, when the $X_i$'s are independent ($\rho = 0$) and when they are weakly dependent ($\rho = 0.10$). Figures 1 and 2 compare the simulated 2-FDRs of these procedures for different values $n_0$ when $\rho = 0$ and 0.10, respectively.

Figures 1 and 2 indicate that while our procedure provides uniformly better control of the $k$-FDR than the generalized Hochberg procedure, which we expected, the difference is, however, quite significant when $n_0$ is neither very small nor very large. Compared to the original BH procedure, we notice that our procedure works quite well when the test statistics are independent or close to being independent and $n_0$ is relatively large, that is, when relatively few of the null hypotheses are expected to be false. With increasing dependence among the test statistics, our procedure loses its edge over the BH procedure for small $n_0$.

**4. Concluding remarks.**   Generalizing traditional error rates to make them more appropriate in situations where one is willing to tolerate more than one false rejection and is wishing to increase the ability of procedures to detect false null hypotheses using error rates that allow such rejections has become an increasingly important idea in multiple testing, because of its relevance in testing a large number of hypotheses, as in microarray studies. The work done in this article makes an important contribution in this area. In addi-tion to generalizing previous work on the $k$-FWER, the notion of FDR has been generalized for the first time and procedures controlling it have been developed in this article. We believe we have opened the door in this article for further research in multiple testing, particularly toward developing the theory and methodology of false discovery rate. There are several interesting generalizations of results related to the original concept of FDR that could potentially be developed along the line of this article.

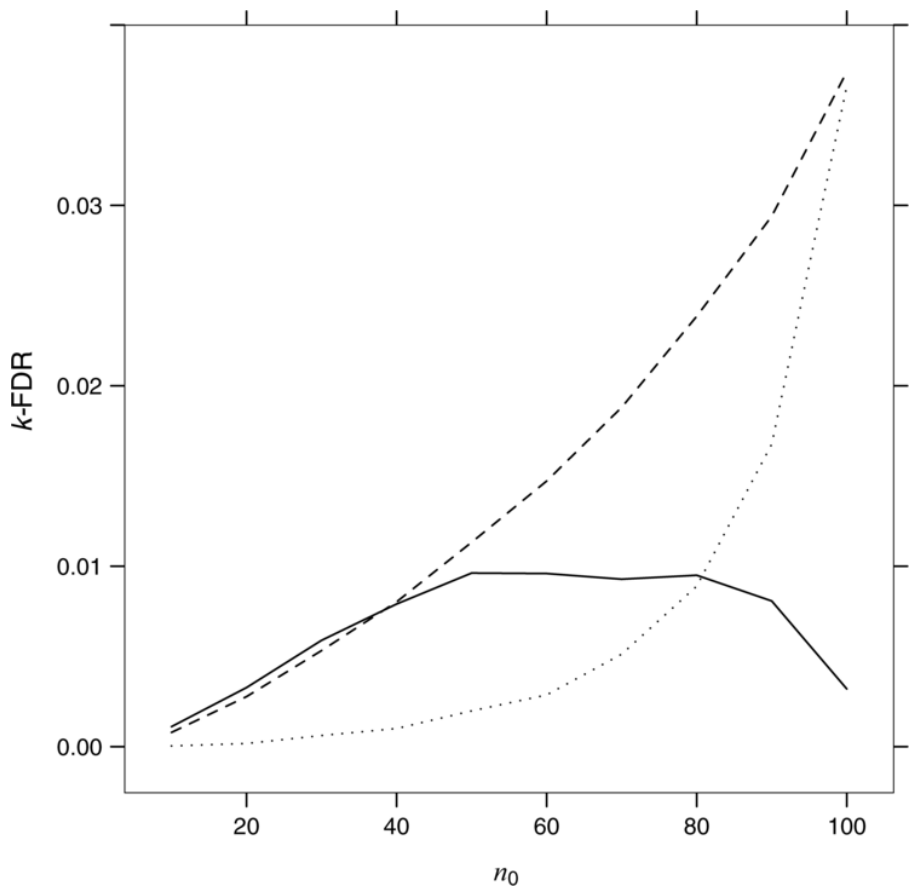

FIG. 1.   *Comparison of 2-FDRs with $\rho = 0.0$ (Generalized BH: - - - ; Generalized Hochberg: $\cdots$; Original BH: —).*

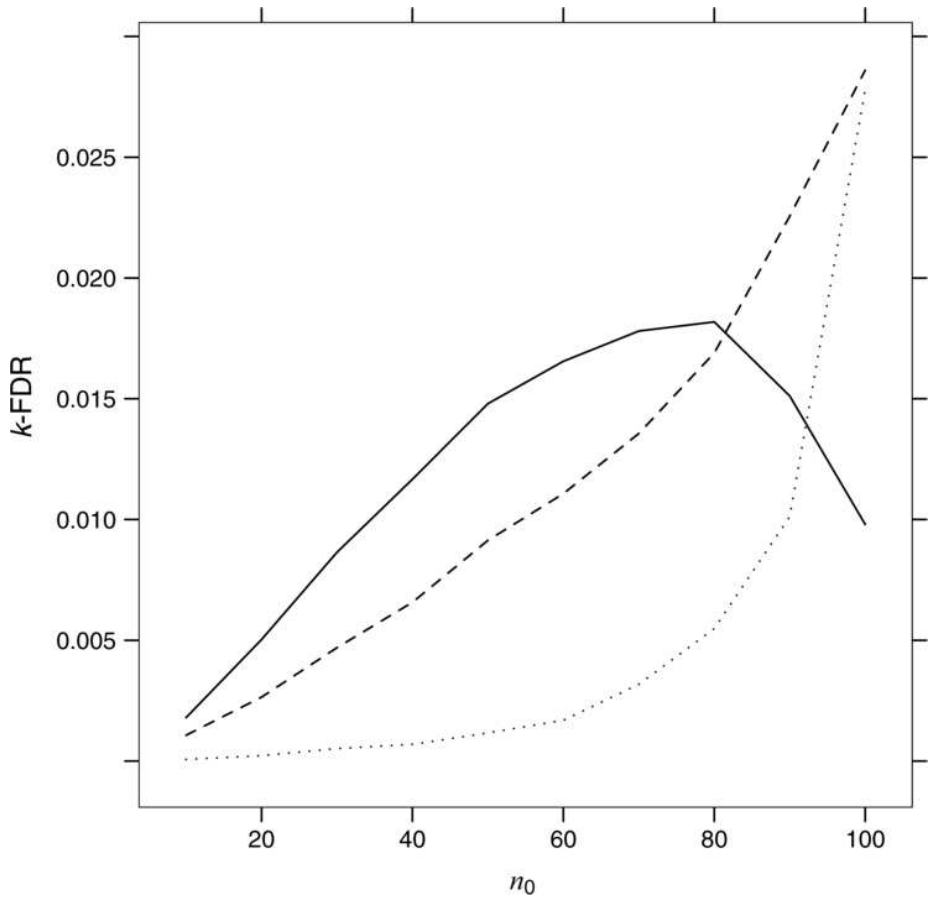

FIG. 2. *Comparison of 2-FDRs with $\rho = 0.10$ (Generalized BH: - - - ; Generalized Hochberg: $\cdots$; Original BH: ---).*

**Acknowledgments.** The author thanks Zijiang Yang for doing the numerical calculations and a referee for giving helpful suggestions.

## REFERENCES

[1] BENJAMINI, Y. and HOCHBERG, Y. (1995). Controlling the false discovery rate: A practical and powerful approach to multiple testing. *J. Roy. Statist. Soc. Ser. B* **57** 289–300. MR1325392

[2] BENJAMINI, Y. and YEKUTIELI, D. (2001). The control of the false discovery rate in multiple testing under dependency. *Ann. Statist.* **29** 1165–1188. MR1869245

[3] BENJAMINI, Y. and YEKUTIELI, D. (2005). False discovery rate-adjusted multiple confidence intervals for selected parameters (with discussion). *J. Amer. Statist. Assoc.* **100** 71–93. MR2156820

[4] FAN, J., HALL, P. and YAO, Q. (2006). To how many simultaneous hypothesis tests can normal, Student's $t$ or bootstrap calibration be applied? Unpublished manuscript.

[5] GENOVESE, C. and WASSERMAN, L. (2002). Operating characteristics and extensions of the false discovery rate procedure. *J. R. Stat. Soc. Ser. B Stat. Methodol.* **64** 499–517. MR1924303

[6] GENOVESE, C. and WASSERMAN, L. (2004). A stochastic process approach to false discovery control. *Ann. Statist.* **32** 1035–1061. MR2065197

[7] HOCHBERG, Y. (1988). A sharper Bonferroni procedure for multiple tests of significance. *Biometrika* **75** 800–802. MR0995126

[8] HOLM, S. (1979). A simple sequentially rejective multiple test procedure. *Scand. J. Statist.* **6** 65–70. MR0538597

[9] KARLIN, S. and RINOTT, Y. (1980). Classes of orderings of measures and related correlation inequalities. I. Multivariate totally positive distributions. *J. Multivariate Anal.* **10** 467–498. MR0599685

[10] KORN, E., TROENDLE, J., MCSHANE, L. and SIMON, R. (2004). Controlling the number of false discoveries: Application to high-dimensional genomic data. *J. Statist. Plann. Inference* **124** 379–398.

[11] LEHMANN, E. L. and ROMANO, J. P. (2005). Generalizations of the familywise error rate. *Ann. Statist.* **33** 1138–1154. MR2195631

[12] MEINSHAUSEN, N. and RICE, J. (2006). Estimating the proportion of false null hypotheses among a large number of independently tested hypotheses. *Ann. Statist.* **34** 373–393. MR2275246

[13] ROMANO, J. P. and SHAIKH, A. M. (2006). Stepup procedures for control of generalizations of the familywise error rate. *Ann. Statist.* **34** 1850–1873. MR2283720

[14] SARKAR, S. K. (1998). Some probability inequalities for ordered MTP$_2$ random variables: A proof of the Simes conjecture. *Ann. Statist.* **26** 494–504. MR1626047

[15] SARKAR, S. K. (2002). Some results on false discovery rate in stepwise multiple testing procedures. *Ann. Statist.* **30** 239–257. MR1892663

[16] SARKAR, S. K. (2004). FDR-controlling stepwise procedures and their false negatives rates. *J. Statist. Plann. Inference* **125** 119–137. MR2086892

[17] SARKAR, S. K. (2006). False discovery and false nondiscovery rates in single-step multiple testing procedures. *Ann. Statist.* **34** 394–415. MR2275247

[18] SARKAR, S. K. (2007). Generalizing Simes' test and Hochberg's stepup procedure. *Ann. Statist.* To appear.

[19] SARKAR, S. K. and CHANG, C.-K. (1997). The Simes method for multiple hypothesis testing with positively dependent test statistics. *J. Amer. Statist. Assoc.* **92** 1601–1608. MR1615269

[20] SIMES, R. J. (1986). An improved Bonferroni procedure for multiple tests of significance. *Biometrika* **73** 751–754. MR0897872

[21] Storey, J. D. (2002). A direct approach to false discovery rates. *J. R. Stat. Soc. Ser. B Stat. Methodol.* **64** 479–498. MR1924302

[22] Storey, J. D. (2003). The positive false discovery rate: A Bayesian interpretation and the $q$-value. *Ann. Statist.* **31** 2013–2035. MR2036398

[23] van der Laan, M., Dudoit, S. and Pollard, K. (2004). Augmentation procedures for control of the generalized family-wise error rate and tail probabilities for the proportion of false positives. *Stat. App. Gen. Mol. Biol.* **3** Article 15. MR2101464

Department of Statistics
Fox School of Business and Management
Temple University
Philadelphia, Pennsylvania 19122
USA
E-mail: sanat@temple.edu